\documentclass[10 pt,a4paper]{article}
\usepackage[toc,page]{appendix}
\usepackage{amsmath}
\usepackage{tocbibind}
\usepackage{stmaryrd}   
\usepackage{amsmath}
\usepackage{verbatim}
\usepackage{mathrsfs}
\usepackage{mathtools}
\usepackage{amsthm} 
\usepackage{amssymb} 
\usepackage{mdframed}
\usepackage[margin=2.5 cm]{geometry}
\usepackage{enumerate}
\usepackage{bm}
\usepackage[colorlinks={true},linkcolor={blue},citecolor={blue}, urlcolor={blue}]{hyperref}
\usepackage[noabbrev,nameinlink]{cleveref}
\usepackage{comment}
\usepackage{tikz}
\usepackage{tikz-cd}
\usetikzlibrary{intersections}
\usetikzlibrary{spy} 
\theoremstyle{plain}

\newtheorem{thm}{Theorem}[section]
\newtheorem{prop}[thm]{Proposition}

\theoremstyle{definition}
\newtheorem{defn}[thm]{Definition}

\newtheorem{expl}[thm]{Example}

\newtheorem{rem}[thm]{Remark}

\newtheoremstyle{named}%
    {}{}{\itshape}{}{\bfseries}{.}{.5em}{\thmnote{#3}}
\theoremstyle{named}

 \title{A note on some Diophantine inequalities over adelic curves}
\author{Paolo Dolce}
\date{}

\newcommand{\Addresses}{{
  \bigskip
  \footnotesize
}
  
  P.~Dolce, \textsc{Ben Gurion University of the Negev and Westlake University}\par\nopagebreak
  \textit{E-mail address}: \texttt{dolce@westlake.edu.cn}}

\begin{document}

\maketitle
\abstract{Without assuming the Northcott property we provide an upper bound on the number of ``big solutions'' of a special system of Diophantine inequalities over proper adelic curves. This system is interesting since it is a stronger version Roth's inequality for adelic curves.}

\makeatletter
\@starttoc{toc}
\makeatother

\section{Introduction}
Several theorems prove the finiteness of the solutions of certain Diophantine equations and inequalities, both over number fields and over function fields. Regarding any Diophantine problem with finitely many solutions we talk about ``effective solution'' when we are able to find an explicit bound on the height of the solutions,  whereas we refer to ``quantitative solution'' when we are able to bound only the number of solutions. In general one can deduce the latter from the former. Diophantine problems are very often much easier on function fields than on number fields, and  a classical example is given by the $S$-unit equation: indeed for function fields we have a simple proof for the effective solution (see for instance \cite[Theorem 3.16]{Z}) whereas for number fields we only have quantitative results (see \cite{eve}) or some partial effective results (see for instance \cite{BG}).   

Roughly speaking, \emph{adelic curves} are fields whose absolute values are parameterized by a measure space; they have been introduced by Chen and Moriwaki in \cite{Ch-Mo}. They are a generalisation of number fields and function fields, but they also generalise a very important and less understood family of fields, i.e. arithmetic function fields. Such fields are finite generated extension of $\mathbb Q$ endowed with a set of absolute values coming from  an arithmetic variety completed in the Arakelov sense (see \cite[Section 3]{vo}). Therefore, it is interesting  to formulate and solve Diophantine problems for adelic curves, because this means having unifying proofs of theorems that belong to apparently different frameworks.

For instance, two generalisations of Roth's theorem for adelic curves have been proved in \cite{dz}. Some extra properties on the adelic curves are needed in order to ensure the validity of Roth's theorem(s), nevertheless the proved results include the classical Roth's theorem for number fields, function fields and also the the more recent version of Roth's theorem for arithmetic function fields of \cite{vo}. At the time of writing we don't know of any quantitative version of Roth's theorem for adelic curves, whereas several quantitative results are known only for number fields and  function fields: \cite{DR}, \cite{bom}, \cite{bom1}, \cite{gr1}, \cite{sc2}, \cite{C}.

 In this short paper we investigate the quantitative solution of certain Diophantine inequalities for adelic curves. Let's first fix some notations that we will keep throughout the whole paper, some more details on the definitions will be given in \Cref{sec:prel}:\\

\textbf{Notations.} $\mathbb K$ is a field of characteristic $0$ and $\mathbb X=(\mathbb K,\Omega,\phi)$ is a proper adelic curve. The adjective ``proper'' means that we have a product formula. $\Omega$ denotes the measure space, with measure $\mu$, that parameterises the absolute values of $\mathbb K$ through the function $\phi$.  The logarithmic height on $\mathbb X$ is denoted by $h=h_{\mathbb X}$, and we will also make use of the exponential height $H=e^h$. We fix  a subset $S\subset\Omega$  of finite measure; $\alpha_1,\ldots, \alpha_n$ are distinct elements of $\mathbb K^\times$; $A>0$ and $\varepsilon>0 $ are real constants; $\theta:S\to\mathbb R_{\ge 0}$ is a measurable function such that $\int_S\theta d\mu\ge 2+\varepsilon$.\\

With the above notations fixed, we study the following system of inequalities in the variable $\beta\in\mathbb K^\times$:
\begin{equation}\label{eq:main_system}
\left|\beta-\alpha_i\right|_{\omega}\le (AH(\beta))^{-\theta(\omega)}\,,\quad \text{ for $i=1,\ldots, n$ and $\omega\in S$}\,. 
\end{equation}
Before the presentation of the results let's explain the importance of this Diophantine problem. Assume that $\beta$ satisfies the inequalities of $(\ref{eq:main_system})$, and choose any map $S\ni\omega\mapsto\alpha_\omega\in\{\alpha_1,\ldots,\alpha_n\}$. Then from $(\ref{eq:main_system})$ we get:
\begin{equation}\label{eq:main_and_roth}
\int_S\frac{\log|\beta-\alpha_\omega|_\omega d\mu}{\log A+h(\beta)}\le -\int_S \theta d\mu\le -2-\varepsilon
\end{equation}
Thus $\beta$ is a solution of  (a modified) Roth's inequality for adelic curves (i.e. a good approximant). Recall that \cite[Theorem A]{dz} shows that if $\mathbb X$ satisfies the strong $\mu$-equicontinuity property, then the solutions of Roth's inequality have bounded heights. Moreover, if $n=1$, $\Omega$ has the counting measure, $\theta=2+\varepsilon$, $A=1$ and $S=\{\omega\}$, the system (\ref{eq:main_system}) is exactly the classical Roth's inequality:
\[
|\beta-\alpha|_\omega\le H(\beta)^{-2-\varepsilon}\,.
\]
\begin{rem}
When $\Omega$ has the counting measure Corvaja shows in \cite{C} that  Roth's inequality can be reduced to solving at most $M$ systems of the type  (\ref{eq:main_system}), where $M$ is a number depending on $\#S$. Therefore in this special case an effective/quantitative solution of (\ref{eq:main_system}) would also give an effective/quantitative solution of Roth's inequality. Corvaja's reduction process is based on the following simple fact: for any $\beta\in\mathbb K^\times$ the functions $\omega\mapsto\min(1,|\beta-\alpha_i|_\omega)$ can be seen as points in the hypercube $\mathcal C=[0,1]^{(\#S)}$, then one can subdivide $\mathcal C$ in smaller hypercubes of the same diameter by using a mesh. Then, the points inside $\mathcal C$ can be approximated with arbitrary precision (eventually after refining the mesh) by the finitely many corners of the smaller hypercubes. The usage of some combinatorics on hypercubes is very common for the effective/quantitative solution of Diophantine problems, see for instance \cite{eve} for another application. These techniques clearly don't make sense when $\Omega$ is a general measure space.
\end{rem}
We prove the following quantitative result that gives an upper bound on the solutions of $(\ref{eq:main_system})$ with big height:
\begin{thm}\label{thm:main_th}
Fix $N=\max\left(\lceil 21^2\log 2n\rceil, \lceil \frac{49\log 2n}{\varepsilon^2}\rceil\right)$ and $A\le 4^{2N!}\prod^n_{i=1} H(\alpha_i)^{\frac{2N!}{n}}$. Let $\mathcal B$ be the set of solutions $\beta\in\mathbb K^\times$  of system (\ref{eq:main_system})  satisfying:
 \[
h(\beta)\ge \max\left(\log A, 2\log A +\frac{4\log A-2h(2)-\log 4}{\varepsilon}\right)\;;
\]
then
\[
\#\mathcal B < (N-1)\left(1+\frac{\log 8nN^2N!}{\log(1+\varepsilon/2)}\right)\,.
\]

\end{thm}
The strength of this result is that it holds for any proper adelic curves without extra assumptions, in particular it is independent from the Northcott property (see \Cref{defn:north}). In other words it can be seen as a consequence of the product formula only. Notice that the bound on the number of solutions depends only on $n$ and $\varepsilon$; moreover when $A\le 1$ then  
\[
\max\left(\log A, 2\log A +\frac{4\log A-2h(2)-\log 4}{\varepsilon}\right)\le 0\,,
\]
so \Cref{thm:main_th} gives a bound for \emph{all the solutions} of the system of inequalities.  In the concrete case  of arithmetic function fields this quantitative result seems completely new  and it could lead to interesting applications.

Our proof is  based on the ideas in \cite{C} and it follows from the combination of two gap principles. In general a gap principle for a Diophantine problem is a result that gives a constraint on the distribution of the heights of the solutions (of the problem). The first gap principle of this paper is very simple and is an application of the integral Liouville inequality (recalled in \Cref{prop:liou}). It gives a lower bound, depending on $\varepsilon$, on the ratio of the heights of the solutions. The second gap principle is more complicated since it follows from \cite[Theorem 3.9]{dz} which in turn relies on a highly nontrivial version of Dyson's lemma for multivariate polynomials (see \cite{EV}). Roughly speaking it gives an opposite constraint to the first gap principle, in fact the second gap principle says that system (\ref{eq:main_system}) cannot have ``too many big solutions'' whose heights are ``spaced enough''.

\paragraph{Acknowledgements.} The author wants to express his gratitude to \emph{Francesco Zucconi} for his helpful comments.

\section{Preliminaries}\label{sec:prel}
Here we collect the main definitions from the theory of adelic curves; for more details the reader can consult \cite{Ch-Mo}, \cite{Ch-Mo2} and \cite{dz}.
We will use the following notations:

$$\log^+ x := \max\{0, \log x\}\,,\quad \log^- x := \min\{0, \log x\}\,;\quad\forall x\in \mathbb R_{>0}$$
\begin{defn}\label{AdCu}
 Let $\mathbb K$ be a field of characteristic $0$,  let $M_{\mathbb K}$ be the set of all absolute values of $\mathbb K$ and let $\Omega=(\Omega,\mathcal A, \mu)$ be a  measure space endowed with  a map 
\begin{eqnarray*}
\phi\colon \Omega&\to& M_{\mathbb K}\\
\omega &\mapsto & |\cdot|_\omega:=\phi(\omega)\,.
\end{eqnarray*} 
 such that for any $a\in \mathbb K^\times$, the real valued function $\omega\mapsto \log |a|_\omega$ lies in  $L^1(\Omega,\mu)$.  The triple $\mathbb X=(\mathbb K,\Omega,\phi)$ is called an \emph{adelic curve}; $\Omega$ and $\phi$ are respectively the \emph{parameter space} and the  \emph{parametrization}.  We denote with $\Omega_{\infty}$ the subset of $\Omega$ made  of all elements $\omega$ such that $|\cdot|_{\omega}$ is archimedean. We set  $\Omega_0:=\Omega\setminus \Omega_{\infty}$.
\end{defn}

It is easy to show that the set $\Omega_{\infty}$ is always measurable \cite[Proposition 3.1.1]{Ch-Mo}.

\begin{defn}
An adelic curve $\mathbb X=(\mathbb K,\Omega,\phi)$ is said to be \emph{proper}  if for any $a\in \mathbb K^\times$:

\begin{equation}\label{eq:pr_for}
\int_{\Omega} \log |a|_\omega\, d\mu(\omega)=0\,.
\end{equation}
\end{defn}

The above property is obviously a generalisation of the product formula. Let's mention some important examples of adelic curves:

\begin{expl}\label{counting}
Any number field with the counting measure on $\Omega$  and  the ``natural'' choice of $\phi$  is a proper adelic curve. Notice that in this case $\phi$ could not be injective, in fact in order to obtain a product formula as in equation (\ref{eq:pr_for}) one usually considers twice the complex embeddings.
\end{expl}

\begin{expl}
An arithmetic function field  $K$ with a big polarisation is a proper  adelic curve. For details see \cite[Section 3]{vo}. 
\end{expl}
\begin{expl} 
A polarised algebraic function field (in $d\ge 1$ variables) over a field of characteristic $0$ can be endowed with a structure of proper adelic curve. For details see  \cite[3.2.4]{Ch-Mo}. 
\end{expl}

\begin{defn}
Let $\mathbb X=(\mathbb K,\Omega,\phi)$ be a proper adelic curve then the \emph{(logarithmic) height} of an element $\beta\in\mathbb K^\times$ is defined as
\[
h(\beta)=h_{\mathbb X}(\beta)=\int_{S}\log^+|\beta|_\omega d\mu(\omega)\,.
\]
We also put $H(\beta):=e^{h(\beta)}$.
\end{defn}
If $\overline{\mathbb K}$ is an algebraic closure of $\mathbb K$ one can produce a natural structure of adelic curve on $\overline{\mathbb K}$, so that the height can be defined for any element algebraic over $\mathbb K$. The details  can be checked on the main references of this section, but in this paper we don't need to deal with the extensions of adelic curves. For $\omega\in\Omega_\infty$, $\vert\cdot\vert_\omega$, by Ostrowski's theorem  there exists a real number $\epsilon(\omega)\in\,]0,1]$ such that $\vert\cdot\vert_\omega=\vert\cdot\vert^{\epsilon(\omega)}$ where on the right we mean the standard euclidean absolute value on $\mathbb C$. Thus, we have a map $\epsilon:\Omega_{\infty}\to ]0,1]$ which can be extended to  $\epsilon:\Omega\to [0,1]$ by putting $\epsilon_{|\Omega_0}:=0$. For instance,  for an archimedean $\vert\cdot\vert_\omega$ we have $\log \vert 2\vert_\omega=\epsilon(\omega)\log 2$, therefore we obtain the explicit expression of the function  $\epsilon$ on the whole $\Omega_\infty$:
$$\epsilon(\omega)=\frac{\log^+\vert2\vert_\omega}{\log 2}\,.$$
Clearly $\epsilon$ is a measurable function. Moreover we always normalise the measure on $\Omega$ in order to satisfy the inequality $h(2)\le \log 2$. One can prove several standard properties for the height on $\mathbb X$, but here we will need only the following one:
\begin{prop}[Integral Liouville inequality]\label{prop:liou}
For any $\alpha,\beta\in\mathbb K^\times$  and any measurable $U\subseteq \Omega$ we have:
\[ 
\int_U\log \vert \alpha-\beta\vert_\omega d\mu(\omega)\ge -\log 2- h(\alpha)-h(\beta)
\]
\end{prop}
\proof
See \cite[Proposition 1.13]{dz}.
\endproof
Finally, let's recall the notion of Northcott property for adelic curves
\begin{defn}\label{defn:north}
A proper adelic curve $\mathbb X=(\mathbb K,\Omega, \mu)$ satisfies the \emph{Northcott property} if for any $C\in\mathbb R$ the set $\{\alpha\in\mathbb K\colon h(a)\}\le C$ is finite.
\end{defn}
We point out once again that this notion won't be needed for our finiteness results.

\section{Quantitative result}

\subsection{First gap principle}
We show that for any  two distinct solutions of (\ref{eq:main_system}) with big heights, the ratio between the heights (if it is $\ge 1$) can be bounded from below by a constant depending only on $\varepsilon$. In other words the heights of such solutions must be separated by ``big enough gaps''.  
\begin{prop}\label{prop:gap1}
Let $\beta_1,\beta_2\in\mathbb K^\times$ two distinct solutions of the system of inequalities (\ref{eq:main_system}) such that:
\begin{equation}\label{eq: cond_h_1}
h(\beta_2)\ge h(\beta_1)\ge 2\log A +\frac{4\log A-2h(2)-\log 4}{\varepsilon}\,,   
\end{equation}
then
\[
\frac{h(\beta_2)}{h(\beta_1)}\ge1+\frac{\varepsilon}{2}\,.
\]    
\end{prop}
\proof
If $\omega$ is archimedean, by the triangle inequality we have
\[
|\beta_1-\beta_2|^{1/\epsilon(\omega)}_\omega\le |\beta_1-\alpha|^{1/\epsilon(\omega)}_\omega+|\beta_2-\alpha|^{1/\epsilon(\omega)}_\omega\le 2(AH(\beta_1))^{-\theta(\omega)/\epsilon(\omega)}
\]
where $\epsilon(\omega)$ is the function defined in the previous section. So
\[
|\beta_1-\beta_2|_\omega\le 2^{\epsilon(\omega)}(AH(\beta_1))^{-\theta(\omega)}\,.
\]
Then we take the logarithm and the integral over $\Omega_\infty$ and we obtain the following inequality:
\begin{equation}\label{eps1}
\int_{\Omega_\infty}-\log|\beta_1-\beta_2|_\omega\ge(\log A+h(\beta_1))\!\int_{\Omega_\infty}\!\!\theta d\mu-\log 2\!\int_{\Omega_\infty}\!\!\epsilon d\mu \ge (\log A+h(\beta_1))\!\int_{\Omega_\infty}\!\!\theta d\mu- h(2)
\end{equation}

If $\omega$ non-archimedean
\[
|\beta_1-\beta_2|_\omega\le \max\{|\beta_1-\alpha|_\omega,|\beta_2-\alpha|_\omega \}\le (AH(\beta_1))^{-\theta(\omega)}\,.
\]
and after integrating over $\Omega_0$ we get
\begin{equation}\label{eps2}
\int_{\Omega_0}-\log|\beta_1-\beta_2|_\omega\ge(\log A+h(\beta_1))\!\int_{\Omega_0}\!\!\theta d\mu
\end{equation}
At this point we sum (\ref{eps1}) and (\ref{eps2}), and we use the property  $\int_{\Omega}\theta d\mu\ge\int_{S}\theta d\mu\ge 2+\varepsilon$:
\[
\int_\Omega-\log|\beta_1-\beta_2|_\omega\ge (\log A+h(\beta_1))\!\int_{\Omega}\!\!\theta d\mu- h(2)\ge  (2+\varepsilon)(\log A+h(\beta_1))-h(2)
\]
\Cref{prop:liou}  now implies
\[
h(\beta_1)+h(\beta_2)+\log 2\ge (2+\varepsilon)(\log A+h(\beta_1))-h(2)
\]
and dividing by $h(\beta_1)$ we get:
\[
\frac{h(\beta_2)}{h(\beta_1)}\ge 1 +\varepsilon+\frac{(2+\varepsilon)\log A-h(2)-\log 2}{h(\beta_1)}\,.
\]
We conclude thanks to the lower bound of inequality (\ref{eq: cond_h_1}).
\endproof
\begin{rem}
Note that if $\log A\le\frac{2h(2)+\log 4}{2\varepsilon +4}$ the condition on the lower bound of $h(\beta_1)$ is empty and \Cref{prop:gap1} gives a gap bound for all the solutions. 
\end{rem}

\subsection{Second gap principle}
For commodity let's introduce the following notations:  
\begin{defn}
  Let  $I$ be a real bounded interval  with $a=\inf I$, $b=\sup I$ and $1<a<b$. The \emph{logarithmic length} of $I$ is defined as the real number $\log\frac{b}{a}$. If $a<1$ one can anyway define the logarithmic length after an adequate translation of  the interval.
\end{defn}

\begin{defn}
    We say that two vectors  $\bm{\alpha}=(\alpha^{(1)},\ldots,\alpha^{(N)})$, $\bm{\beta}=(\beta^{(1)},\ldots,\beta^{(N)})$  are \emph{componentwise different} if $\alpha^{(j)}\neq\beta^{(j)}$ for  $j=1,\ldots, N$.
\end{defn}
We now need to recall some notions and results from \cite{dz}; we also keep the exact same notations of \cite{dz} (even for the indexes) in order to ease the comparison.
Consider a proper  adelic curve $(\mathbb K,\Omega, \phi)$ and a set of vectors $\bm{\alpha}_1,\ldots,\bm{\alpha}_n,\bm{\beta}\in\mathbb K^N$ which are componentwise different. We construct the following  matrices $T:=T(\bm{\alpha}_1,\ldots\bm{\alpha}_n)\in M(n\times N,\mathbb K)$ and $T(\bm{\beta})\in M(n+1\times N,\mathbb K)$:

\[
T:=\begin{pmatrix}
 \alpha^{(1)}_1& \alpha^{(2)}_1  &\ldots &\alpha^{(N)}_1 \\
 \alpha^{(1)}_2& \alpha^{(2)}_2  &\ldots &\alpha^{(N)}_2\\
\vdots & \vdots & \vdots & \vdots\\
\alpha^{(1)}_n& \alpha^{(2)}_n  &\ldots &\alpha^{(N)}_n
\end{pmatrix};\quad
T(\bm{\beta})=\begin{pmatrix}
 \alpha^{(1)}_1& \alpha^{(2)}_1  &\ldots &\alpha^{(N)}_1 \\
 \alpha^{(1)}_2& \alpha^{(2)}_2  &\ldots &\alpha^{(N)}_2\\
\vdots & \vdots  & \vdots  &\vdots\\
\alpha^{(1)}_n& \alpha^{(2)}_n  &\ldots &\alpha^{(N)}_n\\
\beta^{(1)}& \beta^{(2)}  &\ldots &\beta^{(N)}
\end{pmatrix}
\]
We denote by $\bm{\alpha}^{(j)}$, $j=1,\ldots, N$ the columns of $T$, that is:
\begin{equation*}
\bm{\alpha}^{(j)}=\begin{pmatrix}
 \alpha^{(j)}_1 \\
 \alpha^{(j)}_2\\
\vdots \\
\alpha^{(j)}_n
\end{pmatrix}
\end{equation*} 
and by $\bm{\alpha}_{h}$, for  $h\in\{1,2,\dots ,n\}$, the rows:
\begin{equation*}
\bm{\alpha}_{h}=\begin{pmatrix}
 \mathbb \alpha^{(1)}_h, \alpha^{(2)}_h,\ldots, \alpha^{(N)}_h
 \end{pmatrix}
\end{equation*}
Note that we are asking for the matrices $T$ and $T(\bm{\beta})$ to have componentwise different rows.

\begin{defn}\label{boxing}
For the matrix $T(\bm{\beta})$, consider the following real numbers for any $j=1,\ldots, N$: 
\begin{equation}\label{pivot}
\rho_{j}:=4^{2N!}H(\beta^{(j)})\prod_{h=1}^{n}H(\alpha_h^{(j)})^{\frac{2N!}{n}}
\end{equation}
\begin{equation}\label{pivotdue}
\rho'_{j}:=4^{N!}H(\beta^{(j)})\prod_{h=1}^{n}H(\alpha_h^{(j)})^{\frac{N!}{2n}}
\end{equation}
We say that $T(\bm{\beta})$ satisfies the \emph{$h$-gap property} if the following inequality is satisfied:
$$
\frac{{\rm{log}}\, \rho_j }{{\rm{log}}\,\rho'_{j+1}}<\frac{1}{4nN^2N!},\qquad \forall j=1,\ldots, N-1
$$
\end{defn}

\begin{defn}\label{bounding} Fix a measurable set $S=S_1\sqcup\ldots\sqcup S_n$. We say that an integrable function $\theta\colon S\to\mathbb R_{\ge 0}$ is a \emph{column bounding function for} $T(\bm{\beta})$ on $S$ if 
the following inequality holds:
\begin{equation}\label{boundingbis}
-\frac{1}{\log \rho_{j}}{\log}\vert \alpha^{(j)}_h-\beta^{(j)}\vert_\omega\ge\theta (\omega)\qquad \forall j=1,\ldots,N\,,\; \forall h=1,\ldots,n\,,\; \forall\omega\in S_h 
\end{equation}
\end{defn}

\begin{rem}\label{rem_forthelemma}
Note that if $\theta$ is  column bounding then it follows immediately that $\vert \alpha^{(j)}_h-\beta^{(j)}|_{\omega}\le 1$ for any  $j=1,\dots ,N$, $h=1,\dots,n$, $\omega\in S_h$.
\end{rem}

The following nontrivial theorem is the crucial ingredient for the second gap principle. 

\begin{thm}\label{thm:teoremadue}  Let $(\mathbb K,\Omega, \phi)$ be a proper adelic curve. Assume that the matrix $T(\bm{\beta})\in M(n\times N,\mathbb K)$ satisfies the h-gap property (see Definition \ref{boxing}). Moreover fix $N>21^2\log 2n$ and let $\theta\colon S\to\mathbb R_{\ge 0}$ be a column-bounding function for $T(\bm{\beta})$ (see Definition \ref{bounding}). Then
$$
\int_{ S}\theta d \mu<2+\frac{7\sqrt{\log 2n}}{\sqrt{N}}
$$
\end{thm}
\proof
See \cite[Theorem 3.9]{dz}.
\endproof

The following gap principle informally says that for $A$ not  very big we cannot allow too many solutions of (\ref{eq:main_system}) having heights which are at the same time big and very spaced.

\begin{prop}\label{prop:gap2}
Fix  $N=\max\left(\lceil 21^2\log 2n\rceil, \lceil \frac{49\log 2n}{\varepsilon^2}\rceil\right)$ and $A\le 4^{2N!}\prod^n_{i=1} H(\alpha_i)^{\frac{2N!}{n}}$. Then the solutions $\beta\in\mathbb K^\times$ of (\ref{eq:main_system}) such that $H(\beta)\ge A$ are contained in at most $N-1$ intervals of logarithmic length less or equal than $\Gamma=8nN^2N!$
\end{prop}
\proof
Assume by contradiction that there are $N$ solutions $\beta_1,\ldots,\beta_N$ satisfying
\[
\frac{h(\beta_j)}{h(\beta_{j+1})}<\frac{1}{\Gamma}
\]
and consider the matrix
\[
T(\bm\beta)=\begin{pmatrix}
 \alpha_1& \alpha_1  &\ldots &\alpha_1 \\
 \alpha_2& \alpha_2  &\ldots &\alpha_2\\
\vdots & \vdots & \vdots & \vdots\\
\alpha_n& \alpha_n  &\ldots &\alpha_n\\
\beta_1& \beta_2  &\ldots &\beta_N
\end{pmatrix}
\]
Then
\begin{equation}\label{eq:h_gap}
\frac{h(\beta_j)+2N!\log 4+\sum^n_{i=1}\frac{2N!}{n}h(\alpha_i)}{h(\beta_{j+1})+N!\log 4+\sum^n_{i=1}\frac{N!}{n}h(\alpha_i)}\le \frac{h(\beta_j)+\log A}{h(\beta_{j+1})}\le \frac{2h(\beta_j)}{h(\beta_{j+1})}<\frac{2}{\Gamma}=\frac{1}{4nN^2N!}
\end{equation}
which means that $T(\beta)$ satisfies the $h$-gap property. Moreover since any $\beta_j$ is a solution of (\ref{eq:main_system}) we have that $\theta$ is column bounding for $T(\beta)$, indeed:
\[
-\frac{1}{\log\rho_j}\log|\beta-\alpha_i|_\omega\ge -\frac{1}{h(\beta_j)+\log A}\log|\beta-\alpha_i|_\omega \ge\theta(\omega)\,.
\]
Thus we can apply \Cref{thm:teoremadue} and contradict $\int_S\theta d\mu >2+\varepsilon$ (because of the choice of $N$).
\endproof

\subsection{Conclusion of the proof}
Let's  keep the notations fixed in the introduction of the paper. We use the above gap principles to prove the main result.\\

\textbf{Proof of \Cref{thm:main_th}.}
By \Cref{prop:gap1} the height function $h$ is injective on $\mathcal B$, therefore the cardinality of $\mathcal B$ is equal to  the cardinality of $h(\mathcal B)$. Let $x_1<x_2<\ldots <x_r$ be a finite sequence of real numbers such that  $x_j<x^\gamma_{j+1}$ with $\gamma <1$. We want to bound the number of such elements that can be contained in an interval of logarithmic length equal to $\Gamma>1$.  We clearly have 
\[
\frac{1}{\Gamma}\le\frac{\log x_1}{\log x_r}<\gamma^{r-1}
\]
which means $r < 1-\frac{\log\Gamma}{\log\gamma}$. So by using \Cref{prop:gap1} and \Cref{prop:gap2} we conclude that
\[
\# h(\mathcal B)< (N-1)\left(1+\frac{\log 8nN^2N!}{\log(1+\varepsilon/2)}\right)\,.
\]
\hfill{$\square$}
{}\\
\endproof

\bibliographystyle{alpha}
\bibliography{dioph.bib}

\Addresses
\end{document}